\documentclass[smallextended,referee,envcountsect,]{svjour3}
\smartqed
\usepackage{graphicx}
\usepackage{amsmath,amsfonts,amssymb,amscd}
\usepackage{color}
\usepackage{geometry}
\def\R{\mathbb R}
\def\ox{\bar{x}}
\def\ve{\varepsilon}
\def\B{\mathbb B}
\def\dom{\mbox{\rm dom}\,}
\def\st{\stackrel}
\def\ov{\bar{v}}
\def\ra{\rangle}
\def\la{\langle}
\def\O{\Omega}
\def\N{\mathbb{N}}
\def\epi{\mbox{\rm epi}\,}
\def\gph{\mbox{\rm gph}\,}
\def\oR{\Bar{\R}}
\def\kk{\kappa}
\def\al{\alpha}
\def\lm{\lambda}
\def\gg{\gamma}
\begin{document}

\title{Quadratic Growth and Strong Metric Subregularity of the Subdifferential for a Class of Non-prox-regular Functions}

\author{Nguyen Huy  Chieu \and   Nguyen Thi Quynh Trang \and Ha Anh Tuan}

\institute{Nguyen Huy  Chieu, Corresponding author \at
             Department of Mathematics,  Vinh University\\
               Nghe  An, Vietnam\\
              nghuychieu@gmail.com
                        \and
             Nguyen Thi Quynh Trang \at
             Department of Mathematics,  Vinh University\\
               Nghe  An, Vietnam\\
            nqtrang609@gmail.com
\and
               Ha Anh Tuan  \at
              Faculty of Basic Science, Ho Chi Minh city University of Transport\\
               Ho Chi Minh city, Viet Nam \\
             hatuanhuyhoang@gmail.com
}

\date{Received: date / Accepted: date}

\maketitle

\begin{abstract}
This paper mainly studies  the quadratic growth  and the  strong metric subregularity of the subdifferential  of a function  that  can be represented as the  sum of a function twice differentiable in the extended sense and a  subdifferentially continuous, prox-regular, twice epi-differentiable function. For  such a function, which is not necessarily  prox-regular,  it is shown  that the  quadratic growth, the strong metric subregularity of the subdifferential at a local minimizer,   and the positive definiteness of  the subgradient graphical derivative at a  stationary point are equivalent. In addition,    other characterizations of the quadratic growth and  the strong metric subregularity of the subdifferential are also given.
  Besides, properties of functions twice differentiable in the extended sense are examined.
\end{abstract}

\section{Introduction}

Quadratic growth   is an important property of  extended-real-valued functions, which  plays  a central role  in  optimization  \cite{AG08,AG14,BLN18,BS00,CHNT21,DI15,DL18,DMN14,MMS20,MMS21,MS20}. It can be  used  for justifying  the linear  convergence  of various optimization algorithms \cite{BLN18,DL18,OM21} as well as analyzing perturbations of optimization problems. Especially,  for many favorable  classes of functions,  quadratic growth  is  closely related to critical point stability \cite{AG08,AG14,CHNT21,DI15,MMS20,MMS21,MS20}.

For a proper lower semicontinuous convex function,  as shown by Arag\'on Artacho and Geoffroy~\cite{AG08},  the quadratic growth and   the strong metric subregularity of the subdifferential at a local minimizer are equivalent, and they can be characterized by the positive definiteness of  the subgradient graphical derivative at a  stationary point.

  For an arbitrary  proper lower semicontinuous function,  Drusvyatskiy et al.~\cite{DMN14} proved  the validity of   the quadratic growth under   the strong metric subregularity of the subdifferential   at  a local minimizer.  Drusvyatskiy and Ioffe~\cite{DI15} established  that  the converse  holds whenever  the function  under consideration  is semi-algebraic, and it  may fail if the function is not semi-algebraic.
 It is worth noting that the approach of   \cite{DI15} is based on some facts from semi-algebraic geometry, which might not be available for functions that are not semi-algebraic.

  Using tools of  second-order  variational analysis,  Chieu et al. \cite{CHNT21} showed that   for a proper lower semicontinuous function, the positive definiteness of  the subgradient graphical derivative at a  stationary point guarantees  that the point is a local minimizer and   the subdifferential is strongly metrically subregular, which   implies by \cite{DMN14} that the quadratic growth holds.  Furthermore,
 the  quadratic growth, the strong metric subregularity of the subdifferential at a local minimizer, and the positive definiteness of  the subgradient graphical derivative at a  stationary point are equivalent whenever the  function  is either  subdifferentially continuous, prox-regular, and twice epi-differentiable or variationally convex.

More recent developments in this direction can be found in \cite{MMS20,MMS21,MS20,OM21},  where the authors investigated composite models under certain assumptions on the component functions that make  the composite function  subdifferentially continuous, prox-regular, and twice epi-differentiable.

To the best of our knowledge, all known results on the equivalence relationship between the quadratic growth and  the strong metric subregularity of the subdifferential,  except for the one of  Drusvyatskiy and Ioffe \cite{DI15}, are established only for subclasses of the class of   subdifferentially continuous and prox-regular functions.
This observation leads us to  the question if such an  equivalence relationship is valid for functions that are neither  subdifferentially continuous and prox-regular nor semi-algebraic.

In the current work,  we  study  the quadratic growth  and the  strong metric subregularity of the subdifferential  of  functions  that  can be represented as the  sum of an extended twice differentiable function and a  subdifferentially continuous, prox-regular, twice epi-differentiable function.
This big class of functions encompasses  subdifferentially continuous, prox-regular, twice epi-differentiable functions as well as twice differentiable functions.

For  a function from the just mentioned class, which is not necessarily  prox-regular,  it is shown  that   the  quadratic growth,  the strong metric subregularity of the subdifferential at a local minimizer,  and the positive definiteness of  the subgradient graphical derivative at a  stationary point are equivalent. In addition,  other characterizations of  quadratic growth as well as the strong metric subregularity of the subdifferential are also  given.
Besides,  properties of functions that are twice differentiable in the extended sense are examined.

The rest of the paper is organized as follows. Section 2 collects   notions  from variational analysis that are needed  in the sequel. Section 3 investigates  functions that are twice differentiable in the extended sense. The focuses of this section are on sum rules and chain rules for second subderivative, parabolic subderivative, and   subgradient graphical derivative, which are used for proving the main results reported in Section 4. Section 4 is devoted to the study of quadratic growth and strong metric subregularity of the subdifferential. Here we specially pay the attention to the relationship between these two  properties.  Besides, we are also  interested in characterizations of  quadratic growth and strong metric subregularity of the subdifferential via the second subderivative. Section 5 summarizes the main results of the paper and presents some remarks on this research direction.

\section{Preliminaries}

This section recalls some concepts  and their properties from variational analysis \cite{M1,M18,RW98},
 which are needed for  our   analysis. Unless otherwise stated,  $\R^n$ is a Euclidean space with inner product $\langle\cdot,\cdot\rangle$ and norm $\|\cdot\|,$   and $\overline{\R}:=\R\cup\{\infty\}.$  The closed ball with center $\ox$ and radius $\ve>0$ is denoted by $\B_\ve(\ox):=\{x\in \R^n\ |\ \|x-\ox\|\le\ve\}.$

\begin{definition}{\rm (\cite{M1,M18,RW98}).
Let  $f: \R^n\rightarrow \overline{\R}$ and  let $\bar x\in\dom f:=\big\{x\in \R^n|\; f(x)<\infty\big\}$.
  The {\it proximal subdifferential} of $f$ at $\bar x\in \dom f$ is defined by
$$\partial_pf(\bar x):=\left\{v\in \R^n\ |\ \liminf\limits_{x\to \bar x}\frac{f(x)-f(\bar x)-\langle v, x-\bar x\rangle}{\|x-\bar x\|^2}>-\infty\right\}.$$
The {\it regular subdifferential}  (also called {\it Fr\'echet subdifferential})  of $f$ at $\bar x\in \dom f$ is given  by
$$\widehat{\partial}f(\bar x):=\left\{v\in \R^n\ |\ \liminf\limits_{x\to \bar x}\frac{f(x)-f(\bar x)-\langle v, x-\bar x\rangle}{\|x-\bar x\|}\geq 0\right\}.$$
The {\it limiting  subdifferential} (also called {\it Mordukhovich subdifferential})  of $f$ at $\bar x\in \dom f$ is defined  by
$$\partial f(\bar x):=\left\{v\in \R^n\ |\ \exists x_k\st{f}\to \bar x, v_k\to v\  \mbox{with}\ v_k\in \widehat{\partial}f(x_k)\right\}.$$
If $\bar x\not\in\dom f,$ one puts $\partial f(\bar x)=\widehat\partial f(\bar x)=\partial_p f(\bar x):=\emptyset.$
}\end{definition}

\begin{definition}{\rm  (\cite{RW98}).
A function  $f: \R^n \to \overline{\R}$ is  said to be {\it prox-regular} at $\ox\in \dom f$ for $\ov\in \partial f(\ox)$ if there exist $r,\ve>0$ such that for all $x,u\in \B_\ve(\ox)$ with $|f(u)-f(\ox)|<\ve$ we have
\begin{equation}\label{prox}
f(x)\ge f(u)+\la v, x-u\ra-\frac{r}{2}\|x-u\|^2\quad \mbox{for all}\quad v\in \partial f(u)\cap \B_\ve(\ov).
\end{equation}
 Moreover, $f$ is said to be {\em subdifferentially continuous} at $\ox$ for $\ov$  if whenever $(x_k, v_k)\to (\ox,\ov)$ with $v_k\in \partial f(x_k)$, one has $f(x_k)\to f(\ox)$.
}\end{definition}
From \eqref{prox} it follows that  $ \partial f(u)\cap \B_\ve(\ov)\subset \partial_p f(x)$ whenever $\|u-\bar x\|<\ve$ with $|f(u)-f(\ox)|<\ve$.
Furthermore, if $f$ is  subdifferentially continuous at $\ox$ for $\ov$, then  the inequality   ``$|f(u)-f(\ox)|<\ve$'' in the definition of prox-regularity above can be removed.

The following result is a direct consequence of \eqref{prox}, which is very useful for us to verify the prox-regularity in the sequel.
\begin{lemma} {\rm (\cite[Theorem 13.36]{RW98}).}
\label{lemProxR} If $f:\R^{n}\to\overline{\R}$ is prox-regular and subdifferentially continuous at $\ox$ for $\bar v$ then there exist $r, \epsilon>0$ such that
	\begin{equation}\label{lemeq1}
	\la v_{2}-v_{1},x_{2}-x_{1}\ra\geq -r\left\| x_{2}-x_{1}\right\|^{2},
	\end{equation}
	for every $x_{1},x_{2}\in\mathbb{B}_{\epsilon}(\ox), \ v_{1}\in\partial f(x_{1})\cap\mathbb{B}_{\epsilon}(\bar v), \ v_{2}\in\partial f(x_{2})\cap\mathbb{B}_{\epsilon}(\bar v).$
\end{lemma}

\begin{definition}{\rm  (\cite{RW98}).
Given a function $f:\R^n\to \overline{\R}$ with $f(\bar x)\in \R,$ the {\it subderivative}  of $f$ at $\bar x$ is the function $df(\bar x): \R^n\to [-\infty, \infty]$ defined by
$$df(\bar x)(w)=\liminf\limits_{w'\stackrel{t\downarrow0}\rightarrow w}\frac{f(\bar x+t w')-f(\bar x)}{t}\quad \mbox{for all}\ w\in \R^n.$$
 The {\it second subderivative} of $f$  at $\bar x$ for $v\in \R^n$ and $w\in\R^n$  is given by
\begin{equation}\label{ssd}
d^2f(\bar x|v)(w)=\liminf\limits_{ \begin{subarray}\quad\quad t \downarrow 0\\
w'\longrightarrow w\end{subarray}}\Delta_{t}^2f(\bar x,v)(w'),
\end{equation}
where $\Delta_{t}^2f(\bar x,v)(w'):=\frac{f(\bar x+t w')-f(\bar x)-t \langle v, w'\rangle}{ \frac{1}{2}t^2}.$
\par  Function  $f$ is said to be {\it twice epi-differentiable} at $\bar x\in \R^n$ for $v\in \R^n$ if
for every $w\in \R^n$ and choice of $t_k\downarrow0$ there exists $w_k\to w$ such that
$$\Delta_{t_k}^2f(\bar x, v)(w_k)\to d^2f(\bar x|v)(w).$$

}
\end{definition}

It is well-known that  fully amenable functions \cite{RW98},  including the maximum of finitely many $C^2$-functions,  are     subdifferentially continuous prox-regular and twice epi-differentiable lower semicontinuous  proper functions  \cite[Corollary~13.15 and Proposition~13.32]{RW98}.
\begin{definition}{\rm (\cite{RW98}). Let $\Omega$ be a nonempty subset of $\R^n$ and $\bar x\in \R^n.$
\par $(i)$ The (Bouligand-Severi) {\it tangent cone} to  $\Omega$ at $\bar x\in \Omega$  is  given  by
$$T_\Omega(\bar x):=\big\{v\in\R^n | \, \exists  t_k \downarrow 0, \  v_k\rightarrow v\ \mbox{ with }\  \bar x+t_kv_k\in\Omega \  \forall  k\in \mathbb{N}\big\}.$$
If $\bar x\not\in \O$ then one puts $T_\Omega(\bar x):=\emptyset.$
\par $(ii)$   The {\it second-order tangent set} to $\Omega$ at $\bar x$ for $w\in T_{\Omega}(\bar x)$ is defined by
$$T^2_{\Omega}(\bar x,w)=\left\{u\in\R^n\ |\ \exists t_k\downarrow 0, u_k\rightarrow u\  \mbox{with}\  \bar x+t_kw+\frac{1}{2}t_k^2u_k\in \Omega\ \forall  k\in \mathbb{N}\right\}.$$
$\O$ is called  {\it parabolically derivable} at $\bar x$ for $w\in \R^n$ if $T^2_{\O}(\bar x,w)\not=\emptyset,$  and for each $u\in T^2_{\O}(\bar x,w)$ there exist $\ve>0$ and a  mapping  $\xi: [0,\ve]\to \O$ such that $\xi(0)=\bar x,$ $\xi_+'(0)=w$ and  $\xi_+''(0)=u,$
where
$$ \xi_+'(0):=\lim\limits_{t\downarrow0}\frac{\xi(t)-\xi(0)}{t}\quad \mbox{and}\quad \xi_+''(0):=\lim\limits_{t\downarrow0}\frac{\xi(t)-\xi(0)-t \xi_+'(0)}{\frac{1}{2}t^2}.$$
}
\end{definition}

\begin{definition}{\rm  (\cite{RW98}).
The {\it subgradient graphical derivative} of $f$ at $\ox$ for $\ov\in \partial f(\ox)$ is the set-valued mapping $D(\partial f)(\ox|\ov): \R^n\rightrightarrows \R^n$ defined  by
 \begin{equation*}
 D(\partial f)(\ox|\ov)(w):=\big\{z\, |\, (w,z)\in T_{{\rm gph}\,\partial f}(\ox,\ov)\big\}\quad \mbox{for all}\quad w\in \R^n.
 \end{equation*}
}
\end{definition}
If  $f$ is twice epi-differentiable, prox-regular, subdifferentially continuous at $\ox$ for $\ov$, then  it is known from  \cite[Theorem~13.40]{RW98} that
\begin{equation}\label{Dh}
D(\partial f)(\ox|\ov)=\partial h\quad \mbox{with}\quad h=\frac{1}{2}d^2f(\bar x|\ov).
\end{equation}
\begin{definition}{\rm  (\cite{MS20}).
A function $f: \R^n\to \overline{\R}$ is said to be {\it parabolically regular} at $\bar x$ for $\bar v\in \R^n$ if $f(\bar x)\in \R$ and for all $w$ with $d^2f(\bar x,\bar v)(w)<\infty$ there exist $t_k\downarrow 0$ and $w_k\to w$ such that
\begin{equation}\label{MSeq3.1} \lim\limits_{k\to \infty}\Delta_{t_k}^2f(\bar x,\bar v)(w_k)=d^2f(\bar x,\bar v)(w)\quad \mbox{and}\quad
\limsup\limits_{k\to \infty}\frac{\|w_k-w\|}{t_k}<\infty.
\end{equation}
A nonempty set $\Omega\subset \R^n$ is called {\it parabolically regular} at $\bar x$ for $\bar v$ if its indicator function $\delta_\O$ is parabolically regular at $\bar x$ for $\bar v.$
}\end{definition}

\begin{definition}{\rm  (\cite{RW98}).
Let  $f: \R^n\to\overline{\R},$ $\bar x\in \dom f,$ and $w\in  \R^n$ with $df(\bar x)(w)\in \R.$
\par $(i)$ The {\it parabolic subderivative} of $f$ at $\bar x$ for $w$ with respect to $z$ is
$$d^2f(\bar x)(w|z):=\liminf\limits_{z'\stackrel{t\downarrow 0}\rightarrow z}\frac{f(\bar x+tw+\frac{1}{2}t^2z')-f(\bar x)-tdf(\bar x)(w)}{\frac{1}{2}t^2}.$$
\par $(ii)$   $f$ is said to be {\it parabolically epi-differentiable} at $\bar x$ for $w$ if
$$\dom d^2f(\bar x)(w|\cdot)=\{z\in \R^n \ |\ d^2f(\bar x)(w|z)<\infty\}\not=\emptyset,$$
and for every $z\in \R^n$ and every $t_k\downarrow0$ there exists $z_k\to z$ such that
\begin{equation}\label{MSeq3.7}
d^2f(\bar x)(w|z):=\liminf\limits_{k\to \infty}\frac{f(\bar x+t_kw+\frac{1}{2}t_k^2z_k)-f(\bar x)-t_kdf(\bar x)(w)}{\frac{1}{2}t_k^2}.
\end{equation}
}
\end{definition}
As shown by  Mohammadi and Sarabi \cite[Proposition 3.6]{MS20},    a function $f:\R^n\to \overline{\R}$  with $\bar v\in \partial_pf(\bar x)$ is parabolically regular at $\bar x$ for $\bar v$ if and only if
 \begin{equation}\label{MSeq3.11}
 d^2f(\bar x,\bar v)(w)=\inf\limits_{z\in \R^n}\left\{d^2f(\bar x)(w|z)-\langle z,\bar v\rangle\right\}\quad \mbox{for all}\ w\in K_f(\bar x|\bar v),
 \end{equation}
where $K_f(\bar x|\bar v):=\{w\in \R^n\ |\ df(\bar x)(w)=\langle \bar v, w\rangle\}$ is called {\it the critical cone} of $f$ at $(\bar x,\bar v).$
Furthermore,  if $f$ is  parabolically regular at $\bar x$ for $\bar v$ and $w\in \dom d^2f(\bar x,\bar v)$ then there exists $\bar z \in\dom d^2f(\bar x)(w|\cdot)$ such that
\begin{equation}\label{MSeq3.12}
 d^2f(\bar x,\bar v)(w)=d^2f(\bar x)(w|\bar z)-\langle \bar z,\bar v\rangle.
 \end{equation}

\section{Twice differentiability in the extended sense}

The concept of twice differentiability of functions in the extended sense was introduced by Rockafellar and Wets \cite[Definition 13.1]{RW98},  which came from  the desire to develop second order differentiability at $\bar x$ without having to assume the existence of the first partial derivatives at every point in some neighborhood of $\bar x$.

 This section investigates properties of functions that are twice differentiable in the extended sense, with special attention paying to  sum rules and chain rules of the equality form  for second subderivative,  parabolic subderivative, and subgradient graphical derivative.

\begin{definition}{\rm (\cite{RW98}). Let $f:\R^n\to \overline{\R}$ be finite at $\bar x.$  We say that
\par $(i)$ $f$ is  {\it differentiable} (resp., {\it strictly differentiable}) at $\bar x$ if there exists an $(1\times n)$-matrix $\nabla f(\bar x),$ called the Jacobian (matrix) of $f$ at $\bar x,$
 such that
 $$\lim\limits_{x\to \bar x}\frac{f(x)-f(\bar x)-\nabla f(\bar x)(x-\bar x)}{\|x-\bar x\|}=0 \quad (\mbox{resp,} \lim\limits_{x,u\to \bar x}\frac{f(x)-f(u)-\nabla f(\bar x)(x-u)}{\|x-u\|}=0);$$
\par $(ii)$    $f$ is {\it twice differentiable} at $\ox$ (in the classical sense)  if it is differentiable on a neighborhood $U$ of $\ox$ and there exists a $n\times n$ matrix $\nabla^2 f(\bar x),$ called  the  Hessian  (matrix) of $f$ at $\bar x$, such that
$$\lim\limits_{x\st{U}{\rightarrow} \bar x}\frac{\nabla f(x)-\nabla f(\bar x)-H(x-\bar x)}{\|x-\bar x\|}=0;$$
\par $(iii)$ 	$f$ is {\it twice differentiable at $\ox$  in the extended sense} if it is differentiable at $\ox,$ and there exist a $n\times n$ matrix $A,$  a neighborhood $U$ of $\ox$ and a subset $D$ of $U$ with $\mu(U\backslash D)=0$ such that  $f$ is  Lipschitz on $U,$ differentiable at every point in $D,$ and
$$\lim\limits_{x\st{D}{\rightarrow} \bar x}\frac{\nabla f(x)-\nabla f(\bar x)-A(x-\bar x)}{\|x-\bar x\|}=0,$$
where $\mu$ denotes the Lebesgue measure on $\R^n.$
	This matrix $A$, necessarily unique, is then called the Hesian (matrix) of $f$ at $\ox$ in the extended sense and is likewise denoted by $\nabla^{2}f(\ox).$		
}
\end{definition}

It is known \cite[Theorem 13.51]{RW98}  that a $C^2$-lower function (and thus a $C^{1,1}$-function) on an open set $\mathcal O\subset\R^n$ is twice differentiable in the extended sense almost everywhere in $\mathcal{O},$ with  extended Hessian being symmetric where they exist.

It is easy to see  that if $f$ is  twice differentiable at $\ox$ then it is twice differentiable at $\ox$  in the extended sense, and the Hessian and the extended Hessian coincide.

The following example shows that there exists a function that is twice differentiable in the extended sense, but  neither  twice differentiable nor prox-regular.

\begin{example}\label{ex}{\rm (Extended twice differentiability does not imply  either  twice differentiability or prox-regularity). Consider the  function   $g:\R\rightarrow\R$  given by
		$$g(x)=\begin{cases}
		x^{10/3}\cos\frac{1}{x}+x^{4}\quad\ &\mbox{ if } x\geq 1,\\
		x^{10/3}\cos\frac{1}{x}+\frac{(2n+1)(2n^2+2n+1)}{n^3(n+1)^3}x+\frac{1}{(n+1)^3}-\frac{1}{n^3} \quad\ &\mbox{ if } x\in \big[\frac{1}{n+1},\frac{1}{n}\big), \ n=1,2,...\\
		0 \quad\ &\mbox{ if } x=0,\\
		g(-x) \quad\ &\mbox{ if } x<0,
		\end{cases}$$
		
		${\it Claim 1:}$ $g$ is twice differentiable at $\ox=0$ in the extended sense, but it is not twice differentiable at $\ox$ in the classical sense.		
		Indeed, we see that $g$ is differentiable at $\ox$, and
		$$\nabla g(x)=\begin{cases}
		\frac{10}{3}x^{7/3}\cos\frac{1}{x}+x^{4/3}\sin\frac{1}{x}+4x^3 &\mbox{ if } x>1,\\
		\frac{10}{3}x^{7/3}\cos\frac{1}{x}+x^{4/3} \sin\frac{1}{x}+\frac{(2n+1)(2n^2+2n+1)}{n^3(n+1)^3} &\mbox{ if } x\in \big(\frac{1}{n+1},\frac{1}{n}\big), \ n=1,2,...\\
		0 \quad\ &\mbox{ if } x=0,\\
		-\nabla g(-x) \quad\ &\mbox{ if } \ x\in(-\infty,0)\setminus\big\{-\frac{1}{n}| \ n\in\mathbb{N}^{*}\big\}.
		\end{cases}$$
		Put $U=(-1,1),$   $D=(-1,1)\setminus\{\frac{1}{n}| \ n\in\mathbb{Z}^{*}\},$ and $A=0.$  Then  $\mu(U\backslash D)=0,$ $g$ is  Lipschitz on $U$ with constant $\kappa=1,$ and differentiable at every point in  $D,$ where $\mu$ is  the Lebesgue measure on $\R.$ Furthermore, for each $x\in \big(\frac{1}{n+1},\frac{1}{n}\big)$ with $n\in \mathbb N^*$  we have
		$$\begin{array}{rl}\left|\frac{\nabla g(x)-\nabla g(\bar x)-A(x-\bar x)}{|x-\bar x|}\right|&\leq\left| \frac{10}{3}x^{4/3}\cos\frac{1}{x}+x^{1/3}\sin\frac{1}{x}\right| +\frac{(2n+1)(2n^2+2n+1)}{n^3(n+1)^3|x|}\\
		&\leq \frac{10}{3}x^{4/3}  +x^{1/3}  +\frac{(2n+1)(2n^2+2n+1)}{n^3(n+1)^2}\to 0\quad \mbox{as}\ n\to \infty.  \end{array}$$
	Combing this with  $\nabla g(x)=-\nabla g(-x)$ for all  $x\in(-\infty,0)\setminus\big\{-\frac{1}{n}| \ n\in\mathbb{N}^{*}\big\},$ we get
		$$\lim\limits_{x\st{D}{\rightarrow} \bar x}\frac{\nabla g(x)-\nabla g(\bar x)-A(x-\bar x)}{|x-\bar x|}=0.$$
		Hence,   $g$ is twice differentiable at $\ox$ in the extended sense. On the other hand, since $g$ is not differentiable at each point $\frac{1}{n}$ with $n\in \mathbb Z^*,$   $g$ is not twice differentiable at $\bar x$ in the classical sense.  Therefore,  the extended twice differentiability does not imply the classical twice differentiability.
		
		${\it  Claim  2:}$ $g$ is not prox-regular at $\ox$ for $\bar v=0.$ Fix   $r>0$  and put  $u_{k}=\frac{1}{2k\pi}, \ x_{k}=\frac{1}{\frac{\pi}{2}+{2k\pi}}$ for every $k\in \N^*.$  Then for each  $k\in \N^*$  there exist $m_k, n_k\in\mathbb{N}^{*}$ such that $u_{k}\in\big(\frac{1}{m_k+1}, \frac{1}{m_k}\big)$ and $x_{k}\in\big(\frac{1}{n_k+1}, \frac{1}{n_k}\big).$ This implies that  $2k\pi<m_k+1$ for all $k.$   So we have		
		$$\begin{array}{rl}
		&\la \nabla g(u_{k})-\nabla g(x_{k}), u_{k}-x_{k}\ra+r\left| u_{k}-x_{k}\right|^{2}\\
&=\Big(\frac{10}{3}\frac{1}{(2k\pi)^{7/3}}+\frac{(2m_k+1)(2m_k^2+2m_k+1)}{m^3_k(m_k+1)^3}-\frac{1}{\big(\frac{\pi}{2}+2k\pi\big)^{4/3}}
-\frac{(2n_k+1)(2n_k^2+2n_k+1)}{n_k^3(n_k+1)^3}\Big)\big(\frac{1}{2k\pi}-\frac{1}{\frac{\pi}{2}+{2k\pi}}\big)\\
		& \ \ \ \ +r\big(\frac{1}{2k\pi}-\frac{1}{\frac{\pi}{2}+{2k\pi}}\big)^{2}\\		
&\leq\frac{\pi}{4k\pi\big(\frac{\pi}{2}+2k\pi\big)}
\Big(\frac{10}{3}\frac{1}{(2k\pi)^{7/3}}+\frac{5}{(m_k+1)^3}-\frac{1}{\big(\frac{\pi}{2}+2k\pi\big)^{4/3}}
+r\frac{\pi}{4k\pi\big(\frac{\pi}{2}+2k\pi\big)}\Big)\\
&\leq\frac{\pi}{4k\pi\big(\frac{\pi}{2}+2k\pi\big)}
\Big(\frac{10}{3}\frac{1}{(2k\pi)^{7/3}}+\frac{5}{(2k\pi)^3}-\frac{1}{\big(\frac{\pi}{2}+2k\pi\big)^{4/3}}
+r\frac{\pi}{4k\pi\big(\frac{\pi}{2}+2k\pi\big)}\Big)<0,\\		
		\end{array}$$
		for all  $k$ large enough.  Note that $\lim\limits_{k\to\infty}\big(u_k,\nabla g(u_{k})\big)=\lim\limits_{k\to\infty}\big(x_k,\nabla g(x_{k})\big)=(0,0).$  Therefore, by Lemma \ref{lemProxR},  $f$ is not prox-regular at $\ox=0$ for $\bar v=0.$
}\end{example}

\begin{example}\label{ex2}{\rm (Twice differentiability does not imply   prox-regularity).
		Consider the function  $f:\R\to\R$  defined by
	$$f(x):=\int_{0}^{x}g(t)dt\quad \mbox{where}\quad 		g(x)=\begin{cases}
        x^{2}\mbox{sin}\frac{1}{x^{2}} \quad\ &\mbox{ if } x\not=0,\\
		0 \quad\ &\mbox{ if } x=0.
		\end{cases}$$
		We see that $\nabla f(x)=g(x)$ for all $x\in \R,$ and $f$ is twice differentiable at every point in $\R$ with
		$$\nabla^{2}f(x)=\nabla g(x)=\begin{cases}
		2x\sin\frac{1}{x^{2}}-\frac{2}{x}\cos\frac{1}{x^{2}}\quad\ &\mbox{ if } x\not=0,\\
		0 \quad\ &\mbox{ if } x=0.
		\end{cases}$$
  We next  prove  that $f$ is not prox-regular at $\bar x:=0$ for $\nabla f(\bar x)=0.$  Arguing by contradiction, suppose that $f$ is prox-regular at $0$ for $0.$ Then, by Lemma~\ref{lemProxR}, there exist $r, \epsilon>0$ such that
	\begin{equation}\label{lemeq1a}
	\la g(u)-g(x), u-x\ra= \la \nabla f(u)-\nabla f(x), u-x\ra\geq -r|u-x|^{2},
	\end{equation}
	for every $u,x\in\mathbb{B}_{\epsilon}(\ox).$ Thus, for each $k\in \mathbb N^*$ sufficiently large,  choosing $u_{k}=\frac{1}{\sqrt{2k\pi}}$ and  $x_{k}=\frac{1}{\sqrt{\frac{\pi}{2}+{2k\pi}}},$ we have
		$$\begin{array}{rl}
		&\la g(u_{k})-g(x_{k}), u_{k}-x_{k}\ra+r\left| u_{k}-x_{k}\right|^{2}\\
&=-\frac{1}{\frac{\pi}{2}+{2k\pi}}\Big(\frac{1}{\sqrt{2k\pi}}-\frac{1}{\sqrt{\frac{\pi}{2}+{2k\pi}}}\Big)+r\Big(\frac{1}{\sqrt{2k\pi}}-\frac{1}{\sqrt{\frac{\pi}{2}+{2k\pi}}}\Big)^{2}\\
		&=\Big(\frac{1}{\sqrt{2k\pi}}-\frac{1}{\sqrt{\frac{\pi}{2}+{2k\pi}}}\Big)\Bigg(-\frac{1}{\frac{\pi}{2}+{2k\pi}}+r\Big(\frac{1}{\sqrt{2k\pi}}-\frac{1}{\sqrt{\frac{\pi}{2}+{2k\pi}}}\Big)\Bigg)\\
		&=\Big(\frac{1}{\sqrt{2k\pi}}-\frac{1}{\sqrt{\frac{\pi}{2}+{2k\pi}}}\Big)
\Bigg(-\frac{1}{\frac{\pi}{2}+{2k\pi}}+r\frac{\sqrt{\frac{\pi}{2}+{2k\pi}}-\sqrt{2k\pi}}{\sqrt{2k\pi}\sqrt{\frac{\pi}{2}+{2k\pi}}} \Bigg)\\
		&=\Big(\frac{1}{\sqrt{2k\pi}}-\frac{1}{\sqrt{\frac{\pi}{2}+{2k\pi}}}\Big)\Bigg(-\frac{1}{\frac{\pi}{2}+{2k\pi}}+\frac{\pi r}{2\sqrt{2k\pi}\sqrt{\frac{\pi}{2}+{2k\pi}}\Big(\sqrt{\frac{\pi}{2}+{2k\pi}}+\sqrt{2k\pi}\Big)} \Bigg)<0.
		\end{array}$$		
	 This contradicts   \eqref{lemeq1a} since $\lim\limits_{k\to \infty}u_k=\lim\limits_{k\to \infty}u_k=0=\bar x$. Therefore, $f$ is not prox-regular at $\ox=0$ for $\bar v=0.$
	}
\end{example}

The following lemma collects some properties of extended twice differentiable functions that will be used in the sequel.
\begin{lemma}{ \rm (\cite[Theorem~13.2]{RW98}).}\label{lm1}
	Let  $f:\R^{n}\to\overline \R$ be  twice differentiable at a point $\ox$ in the extended sense.  Then $\partial f(\ox)=\{\nabla f(\bar x)\}$ and there exists a neigborhood $U$ of $\bar x$ such that
	\begin{equation}\label{SD}
	\emptyset\not=\partial f(x)\subset \nabla f(\bar x) + \nabla^2f(x-\ox) +o(\|x-\ox\|)\mathbb{B},
	\end{equation}
	for every $x\in U.$    Furthermore,  $f$ is  strictly  differentiable at $\ox,$ and
\begin{equation}\label{QE}f(x)=f(\bar x)+\langle \nabla f(\bar x), x-\bar x\rangle +\frac{1}{2}\langle x-\bar x, \nabla^2f(\bar x)(x-\bar x)\rangle+o(\|x-\bar x\|^2).\end{equation}
 Here $o(t)$ stands for  some function of $t$ with $\lim\limits_{t\to 0}\frac{o(t)}{t}=0.$
\end{lemma}

Naturally, we say a mapping $F:\R^n\to \R^m,$ $x\mapsto \big(F_1(x), F_2(x),...,F_m(x)\big)$ is {\it twice differentiable at $\ox$  in the extended sense} if  $F_k$ is twice differentiable at $\ox$  in the extended sense for every $k=1,2,...,m.$  In the sequel, for such a mapping $F,$  the symbol $\nabla^2F(\bar x)(w,v)$ stands for
$\big(\langle \nabla^2F_1(\bar x)w,v\rangle, \langle \nabla^2F_2(\bar x)w,v\rangle,...,\langle \nabla^2F_m(\bar x)w,v\rangle\big)$ for all $v, w\in \R^n.$

The following theorem  provides  sum rules of equality form for gradient graphical derivative, second subderivative and parabolic subderivative.
\begin{theorem}\label{thm1} Suppose that  $\varphi:\R^{n}\to\overline\R$  is  twice differentiable at $\ox$ in the extended sense,  $\psi:\R^{n}\to\overline\R$  is  proper lower semicontinuous around $\bar x$, and  $\bar v\in \partial(\varphi+\psi)(\bar x)$. Then one has
	\begin{equation}\label{thm1eq1}
	D\partial(\varphi+\psi)(\ox| \bar v)(w)=\nabla^{2}\varphi(\ox)(w)+D\partial \psi\big(\ox| \bar v-\nabla \varphi(\ox)\big)(w),
  	\end{equation}
	\begin{equation}\label{thm1eq2}d^{2}(\varphi+\psi)\big(\ox| \bar{v}\big)(w)=\big\la w,\nabla^{2}\varphi(\ox)w\big\ra+d^{2}\psi\big(\ox| \bar{v}-\nabla \varphi(\ox)\big)(w),\end{equation}
and	\begin{equation}\label{thm1eq3a}d^{2}(\varphi+\psi)(\ox)(w|z)=\big\la w,\nabla^{2}\varphi(\ox)w\big\ra+ \nabla \varphi(\bar x)z+d^{2}\psi(\ox)(w|z),\end{equation}
for every $w\in\R^{n}$ and $z\in \R^n.$
\end{theorem}
\noindent{\bf Proof.} We first prove \eqref{thm1eq1}.   To this end, take any  $w\in \R^n$ and   $z\in	D\partial(\varphi+\psi)(\ox| \bar v)(w).$ Then  there exist sequences $t_{k}\downarrow0$ and $(w_{k},z_{k})\to (w,z)$ such that
$$\bar v+t_{k}z_{k}\in\partial(\varphi+\psi)(\ox+t_{k}w_{k}) \quad  \mbox{for all} \ k\in \mathbb N^*.$$
Since  $\varphi$  is  twice differentiable at $\ox$ in the extended sense, it is  Lipschitz continuous around $\ox,$ and by the sum rule of  subdifferential \cite[Theorem 2.19]{M18} and  \eqref{SD},   we get
$$
\begin{array}{rl}
\partial(\varphi+\psi)(\ox+t_{k}w_{k})&\subset\partial\varphi(\ox+t_{k}w_{k})+\partial \psi(\ox+t_{k}w_{k})\\
&\subset\nabla\varphi(\ox)+t_{k}\nabla^{2}\varphi(\ox)(w_{k})+o(\left\| t_{k}w_{k}\right\|)\mathbb{B}+\partial \psi(\ox+t_{k}w_{k}),
\end{array}
$$
for all $k\in \mathbb N^*$ sufficiently large.
Thus, for such  numbers $k$,  it holds that
$$\big(\bar{v}-\nabla\varphi(\ox)\big)+t_{k}\Big(z_{k}-\nabla^{2}\varphi(\ox)(w_{k})+\frac{o(\left\| t_{k}w_{k}\right\|)}{t_{k}}\Big)\in\partial \psi(\ox+t_{k}w_{k}),$$
or equivalently,
$$\big(\ox,\bar{v}-\nabla\varphi(\ox)\big)+t_{k}\Big(w_{k}, z_{k}-\nabla^{2}\varphi(\ox)(w_{k})+\frac{o(\left\| t_{k}w_{k}\right\|)}{t_{k}}\Big)\in\mbox{gph}\partial\psi.$$
On the other hand,
$$\Big(w_{k}, z_{k}-\nabla^{2}\varphi(\ox)(w_{k})+\frac{o(\left\| t_{k}w_{k}\right\|)}{t_{k}}\Big)\to\big(w,z-\nabla^{2}\varphi(\ox)(w)\big)\quad \mbox{as}\ k\to \infty.$$
Therefore,
$$\big(w,z-\nabla^{2}\varphi(\ox)(w)\big)\in T_{\mbox{gph}\partial\psi}\big(\ox,\bar{v}-\nabla\varphi(\ox)\big).$$
In other words,
$$z-\nabla^{2}\varphi(\ox)(w)\in D\partial\psi(\ox|\bar{v}-\nabla \varphi(\ox)(w).$$
This shows that
\begin{equation}\label{thm1eq3}
D\partial(\varphi+\psi)(\ox| \bar v)(w)\subset\nabla^{2}\varphi(\ox)(w)+D\partial \psi\big(\ox| \bar v-\nabla\varphi(\ox)\big)(w).
\end{equation}
Conversely, by using \eqref{thm1eq3} and noting that $-\varphi$  is also  twice differentiable at $\ox$ in the extended sense with $\nabla^{2}(-\varphi)(\ox)=-\nabla^{2}\varphi(\ox),$
 we have
$$
\begin{array}{rl}
D\partial\psi(\ox| \bar v-\nabla\varphi(\ox))(w)&=D\partial \big(\varphi+\psi+(-\varphi)\big)(\ox| \bar v-\nabla \varphi(\ox))(w)\\
&\subset D\partial \big(\varphi+\psi\big)(\ox| \bar v)(w)+\nabla^{2}(-\varphi)(\ox)(w)\\
&=D\partial \big(\varphi+\psi\big)(\ox| \bar v)(w)-\nabla^{2}\varphi(\ox)(w).
\end{array}
$$
This infers that
\begin{equation}\label{thm1eq4}
\nabla^{2}\varphi(\ox)(w)+D\partial\psi(\ox| \bar v-\nabla \varphi(\ox))(w)\subset D\partial(\varphi+\psi)(\ox| \bar v)(w).
\end{equation}
From \eqref{thm1eq3} and \eqref{thm1eq4} it follows that
$$D\partial(\varphi+\psi)(\ox| \bar v)(w)=\nabla^{2}\varphi(\ox)(w)+D\partial\psi(\ox| \bar v-\nabla \varphi(\ox))(w) \ \mbox{for every} \ w\in\R^{n}.$$

We next justify the validity of  \eqref{thm1eq2}. Take  any $w\in \R^n.$ Since $\varphi$ is twice differentiable at $\bar x$ in the extended sense, by \eqref{QE},  we see that
 $$\big\la w,\nabla^{2}\varphi(\ox)w\big\ra =\lim\limits_{w'\st{t\downarrow 0}{\rightarrow}w}\Delta^{2}_{t}\varphi\big(\ox|\nabla\varphi(\ox)\big)(w').$$
Therefore,
$$\begin{array}{rl}d^{2}(\varphi+\psi)(\ox| \bar{v})(w)&= \liminf\limits_{w'\st{t\downarrow 0}{\rightarrow}w}\Delta^{2}_{t}(\varphi+\psi)(\ox| \bar{v})(w')\\
&=\liminf\limits_{w'\st{t\downarrow 0}{\rightarrow}w}\Big[\Delta^{2}_{t}\big(\ox|\nabla\varphi(\ox)\big)(w')+\Delta^{2}_{t}\psi\big(\ox| \bar v-\nabla \varphi(\ox)\big)(w')\Big]\\
&=\big\la w,\nabla^{2}\varphi(\ox)w\big\ra+\liminf\limits_{w'\st{t\downarrow 0}{\rightarrow}w}\Delta^{2}_{t}\psi\big(\ox| \bar v-\nabla \varphi(\ox)\big)(w')\\
&= \big\la w,\nabla^{2}\varphi(\ox)w\ra+d^{2}\psi\big(\ox|\bar v-\nabla\varphi(\ox)\big)(w).
\end{array}$$

Finally, we show that  \eqref{thm1eq3a} holds. The differentiability of  $\varphi$ at $\bar x$ gives us that
$$\begin{array}{rl}d(\varphi+\psi)(\ox)(w)&=\liminf\limits_{w'\st{t\downarrow 0}\rightarrow w}\frac{(\varphi+\psi)(\bar x+tw')-(\varphi+\psi)(\bar x)}{t}\\
&=\liminf\limits_{w'\st{t\downarrow 0}\rightarrow w}\left[\frac{\varphi(\bar x+tw')-\varphi(\bar x)}{t}+\frac{\psi(\bar x+tw')-\psi(\bar x)}{t}\right]\\
 &=\nabla \varphi(\bar x)w+  \liminf\limits_{w'\st{t\downarrow 0}\rightarrow w}\frac{\psi(\bar x+tw')-\psi(\bar x)}{t}\\
 &= \nabla \varphi(\bar x)w+  d\psi(\ox)(w)\quad \forall w\in \R^n. \end{array}$$
 Since $\varphi$ is twice differentiable at $\bar x$ in the extended sense,  by \eqref{QE}, we get
  $$\lim\limits_{z'\stackrel{t\downarrow 0}\rightarrow z}\frac{\varphi(\bar x+tw+\frac{1}{2}t^2z')-\varphi(\bar x)-t \nabla \varphi(\bar x)w}{\frac{1}{2}t^2}=\big\la w,\nabla^{2}\varphi(\ox)w\big\ra+\nabla \varphi(\bar x)z \quad \forall w\in \R^n, z\in \R^n.$$
  Therefore,
 \begin{equation*}
 \begin{array}{rl}
&d^{2}(\varphi+\psi)\big(\ox\big)(w|z)=\liminf\limits_{z'\stackrel{t\downarrow 0}\rightarrow z}\frac{(\varphi+\psi)(\bar x+tw+\frac{1}{2}t^2z')-(\varphi+\psi)(\bar x)-td(\varphi+\psi)(\bar x)(w)}{\frac{1}{2}t^2}\\
&=\liminf\limits_{z'\stackrel{t\downarrow 0}\rightarrow z}\Big[\frac{\varphi(\bar x+tw+\frac{1}{2}t^2z')-\varphi(\bar x)-t \nabla \varphi(\bar x)w}{\frac{1}{2}t^2}+\frac{\psi(\bar x+tw+\frac{1}{2}t^2z')-\psi(\bar x)-td\psi(\bar x)(w)}{\frac{1}{2}t^2}\Big]\\
&=\big\la w,\nabla^{2}\varphi(\ox)w\big\ra+ \nabla \varphi(\bar x)z +\liminf\limits_{z'\stackrel{t\downarrow 0}\rightarrow z}\frac{\psi(\bar x+tw+\frac{1}{2}t^2z')-\psi(\bar x)-td\psi(\bar x)(w)}{\frac{1}{2}t^2}\\
&=\big\la w,\nabla^{2}\varphi(\ox)w\big\ra+ \nabla \varphi(\bar x)z+d^{2}\psi(\ox)(w|z)\quad \forall w\in \R^n, z\in \R^n.
\end{array}
\end{equation*}
This finishes the proof. $\hfill\Box$

Let $\psi:\R^n\to \overline{\R}$ be finite at $\bar x\in \R^n.$ Assume that there exists a neighborhood $\mathcal{O}$ of $\bar x$  on which $\psi$ can be represented as \begin{equation}\label{eq4.3}\psi(x)=g\circ F(x)\quad \mbox{for all}\  x\in \mathcal{O},\end{equation}
where  $F: \R^n\to \R^m$ is twice differentiable at $\bar x$ in the extended sense,  and $g: \R^m\to \overline{\R}$ is  proper lower semicontinuous, convex,  and Lipschitz continuous around $F(\bar x)$ relative to its domain with constant $\ell\in \R_+,$ that is, there exists a neighborhood $\mathcal{V}$ of $F(\bar x)$ such that
$|g(y_1)-g(y_2)|\leq \ell\|y_1-y_2\|$ for all $y_1,y_2\in \dom g\cap \mathcal{V}.$

  We see that
\begin{equation}\label{MSeq4.4MSQC}
(\dom \psi )\cap\mathcal{O}=\{x\in \mathcal{O}\ |\  F(x)\in \dom g\}.
\end{equation}

Following \cite[Definition 3.2]{MMS20},  the composition $\psi=g\circ F$ is said to satisfy  the {\it metric subregularity qualification condition} (MSQC)  at $\bar x\in \dom \psi$ if there exist a constant $\kappa\in \R_+$ and a neighborhood $U$ of $\bar x$ such that
\begin{equation}\label{MSeq4.5a}
d(x, \dom\psi)\leq \kappa d\big(F(x),\dom g\big)\quad \mbox{for all}\ x\in U.
\end{equation}

\begin{proposition}\label{prop3.6} Let $\psi:\R^n\to \overline{\R}$ be a function that is  represented as  \eqref{MSeq4.4MSQC} with the composition $g\circ F$ satisfying MSQC  at $\bar x.$  Then we have
$$d\psi(\bar x)(w)=dg\big(F(\bar x)\big)\big(\nabla F(\bar x)w\big)\ \mbox{for all}\  w\in \R^n, \  \partial_p\psi(\bar x)=\partial \psi(\bar x)=\nabla F(\bar x)^*\partial g\big(F(\bar x)\big),$$ and
$T_{\dom \psi}(\bar x)=\big\{w\in \R^n\ |\ \nabla F(\bar x)w\in T_{\dom g}\big(F(\bar x)\big)\big\}.$
\par If assume further that $w\in T_{\dom \psi}(\bar x)$ and $g$ is parabolically epi-differentiable at $F(\bar x)$ for $\nabla F(\bar x)w$, then the following assertions hold:
\par $(i)$ $z\in T^2_{\dom \psi}(\bar x, w)\Leftrightarrow \nabla F(\bar x)z+\nabla^2F(\bar x)(w,w)\in T^2_{\dom g}\big(F(\bar x), \nabla F(\bar x)w\big),$
and $\dom\psi$ is parabolically derivable at $\bar x$ for $w;$
\par $(ii)$ $d^2\psi(\bar x)(w|z)=d^2g\big(F(\bar x)\big)\big(\nabla F(\bar x)w|\nabla F(\bar x)z+\nabla^2F(\bar x)(w,w)$ for all $z\in \R^n$;
\par $(iii)$ $\dom d^2\psi(\bar x)(w|\cdot)=T^2_{\dom\psi}(\bar x,w)$;
\par $(iv)$ $\psi$ is parabolically epi-differentiable at $\bar x$ for $w.$
\end{proposition}
\noindent{\it Proof.} Since $F: \R^n\to \R^m$ is twice differentiable at $\bar x$ in the extended sense, by Lemma \ref{lm1}, we get
\begin{equation}\label{1stStb}F(x)=F(\bar x)+\langle \nabla F(\bar x), x-\bar x\rangle +\frac{1}{2} \nabla^2 F(\bar x)(x-\bar x, x-\bar x)+o(\|x-\bar x\|^2),\end{equation}
and   $f$  is strictly differentiable at $\bar x.$
The latter along with the composition $g\circ F$ satisfying  MSQC  at $\bar x$ implies  by \cite[Theorem 3.4]{MMS20} that
\begin{equation}\label{1stSta}d\psi(\bar x)(w)=dg\big(F(\bar x)\big)\big(\nabla F(\bar x)w\big)\ \mbox{for all}\  w\in \R^n,\end{equation}
and by \cite[Theorem 3.6]{MMS20} that
\begin{equation}\label{1stSa}\partial_p \psi(\bar x)\subset \partial \psi(\bar x)=\nabla F(\bar x)^*\partial g\big(F(\bar x)\big).\end{equation}
We next prove that $\nabla F(\bar x)^*\partial g\big(F(\bar x)\big)\subset \partial_p \psi(\bar x).$ To this end, take any $y\in \partial g\big(F(\bar x)\big).$ Since $g$ is convex,  we have $y\in \partial_p g\big(F(\bar x)\big).$  Hence,
$$\begin{array}{rl}&\liminf\limits_{x\to\bar x}\frac{\psi(x)- \psi(\bar x)-\big\langle \nabla F(\bar x)^*y, x-\bar x\big\rangle}{\|x-\bar x\|^2}\\
&=\liminf\limits_{x\to\bar x}\frac{g\big(F(\bar x)+\nabla F(\bar x)(x-\bar x)+\frac{1}{2}\nabla^2 F(\bar x)(x-\bar x,x-\bar x)+o(\|x-\bar x\|^2)\big)- g\big(F(\bar x)\big)-\big\langle y, \nabla F(\bar x)(x-\bar x)\big\rangle}{\|x-\bar x\|^2}\\
&= \liminf\limits_{x\to\bar x}\frac{g\big(F(\bar x)+\Delta(x)\big)-g\big(F(\bar x)\big)-\big\langle y, \Delta(x)\big\rangle+\big\langle y, \frac{1}{2}\nabla^2 F(\bar x)(x-\bar x,x-\bar x)+o(\|x-\bar x\|^2)\big\rangle}{\|x-\bar x\|^2}\\
&\geq \liminf\limits_{x\to\bar x}\frac{g\big(F(\bar x)+\Delta(x)\big)-g\big(F(\bar x)\big)-\big\langle y, \Delta(x)\big\rangle}{\|x-\bar x\|^2}- \frac{1}{2}\|y\|\cdot\|\nabla^2 F(\bar x)\|>-\infty,\end{array}$$
where $\Delta(x):=\nabla F(\bar x)(x-\bar x)+\frac{1}{2}\nabla^2 F(\bar x)(x-\bar x,x-\bar x)+o(\|x-\bar x\|^2)\rightarrow0$ as $x\to \bar x.$ This shows that $\nabla F(\bar x)^*y\in\partial_p\psi(\bar x)$, and thus
\begin{equation}\label{1stSb}\nabla F(\bar x)^*\partial g\big(F(\bar x)\big)\subset\partial_p\psi(\bar x).\end{equation}
From \eqref{1stSa} and \eqref{1stSb} it follows that
\begin{equation}\label{1stS}\partial_p \psi(\bar x)=\partial \psi(\bar x)=\nabla F(\bar x)^*\partial g\big(F(\bar x)\big).\end{equation}
Furthermore, by \eqref{1stSta} and  \cite[Proposition 2.2]{MS20}, we get
$$\begin{array}{rl} T_{\dom \psi}(\bar x)&= \dom d \psi(\bar x)\\
&=\big\{w\in \R^n\ | \  \nabla F(\bar x)w\in \dom dg\big(F(\bar x)\big)\big\}\\
&=\big\{w\in \R^n\ |\ \nabla F(\bar x)w\in T_{\dom g}\big(F(\bar x)\big)\big\}.\end{array}$$

Let us now suppose further that $w\in T_{\dom \psi}(\bar x)$ and $g$ is parabolically epi-differentiable at $F(\bar x)$ for $\nabla F(\bar x)w.$
Since $g$ is Lipschitz continuous around $F(\bar x)$ relative to its domain, and $\nabla F(\bar x)w\in  T_{\dom g}\big(F(\bar x)\big),$ by \cite[Proposition 4.1]{MS20},  $\dom g$ is parabolically derivable at $F(\bar x)$ for $\nabla F(\bar x)w$. Hence,  the proofs of $(i)$  and  $(ii)-(iv)$  can be, respectively, done  as the ones of \cite[Theorem 4.5]{MMS21} and \cite[Theorem 4.4]{MS20},
where  $F$ was assumed to be  twice differentiable at $\bar x,$ but they actually needed   the quadratic expansion of \eqref{1stStb} and the strict differentiability of $F$ at $\bar x,$ which are valid under the twice differentiability  in the extended sense.
$\hfill\Box$

In order to prove the next proposition  we need  the following lemma whose proof is  the one of   \cite[Proposition 4.6]{MS20}.  For the sake of  completeness we provide the proof with more details.
\begin{lemma} \label{Prop4.6} Suppose  $f: \R^n\to \overline{\R}$ is a proper  lower semicontinuous convex function with
$f(\bar x)\in \R$, $\bar v\in \partial f(\bar x),$ $w\in K_f(\bar x,\bar v),$ and  $f$ is parabolically epi-differentiable at $\bar x$ for $w.$  Then $d^2f(\bar x)(w|\cdot)$ is   proper lower semicontinuous and convex. Furthermore,  $d^2f(\bar x)(w|\cdot)^*(v)=\infty$ whenever  $v\in \R^n\backslash\mathcal{A}(\bar x,w),$ and
$ d^2f(\bar x)(w|\cdot)^*(v)=-d^2f(\bar x,v)(w)$   if  $v\in \mathcal{A}(\bar x,w)$ and  $f$ is  parabolically regular at $\bar x$ for  $v,$
   where  $\mathcal{A}(\bar x,w):=\{v\in \partial f(\bar x)\ |\ df(\bar x)(w)=\langle v, w\rangle\}$ and $d^2f(\bar x)(w|\cdot)^*$ is the Fenchel conjugate of  $d^2f(\bar x)(w|\cdot).$
\end{lemma}
\noindent{\it Proof.} Since $f$ is a lower semicontinuous function, $f(\bar x)\in \R$ and $df(\bar x)(w)=\langle \bar v, w\rangle\in \R,$ by \cite[Proposition 13.64]{RW98}, $d^2f(\bar x)(w|\cdot)$ is lower semicontinuous and
\begin{equation}\label{eqSPsub}d^2f(\bar x)(w|z)-\langle \bar v, z\rangle \geq d^2f(\bar x,\bar v)(w)\quad \forall z\in \R^n.\end{equation}
Noting that $f$ is convex and $\bar v\in \partial f(\bar x)$, we have
\begin{equation}\label{eqSsub} d^2f(\bar x,\bar v)(w)=\liminf\limits_{w'\to w}\frac{f(\bar x+tw')-f(\bar x)-t\langle \bar v, w'\rangle}{\frac{1}{2}t^2}\geq 0.\end{equation}
Thus $d^2f(\bar x)(w|z)>-\infty$ for all $z\in \R^n.$ Combining this with $\dom d^2f(\bar x)(w|\cdot)\not=\emptyset$ (due to the parabolic epi-differentiability of $f$ at $\bar x$ for $w$),  we see that $d^2f(\bar x)(w|\cdot)$ is  a  proper function. By \cite[Example 13.62]{RW98},
$$\epi d^2f(\bar x)(w|\cdot)=T^2_{\epi f}\Big( \big(\bar x, f(\bar x)\big), \big(w, df(\bar x)(w)\big)\Big),$$
and since $f$ is parabolically epi-differentiable at $\bar x$ for $w,$  ${\epi f}$  is parabolically derivable  at $\big(\bar x, f(\bar x)\big)$ for $\big(w, df(\bar x)(w)\big).$  This implies that  $d^2f(\bar x)(w|\cdot)$ is convex since  $f$ is convex.

Take any $v\in \R^n.$ Let us consider the following two cases.

{\it Case 1.} $v\in \mathcal{A}(\bar x,w).$ Then $w\in K_f(\bar x,v)$ and  by  \cite[Proposition 3.6]{MS20}, we have
$$-d^2f(\bar x,v)(w)=-\inf\limits_{z\in \R^n}\{d^2f(\bar x)(w|z)-\langle  v, z\rangle\}=d^2f(\bar x)(w|\cdot)^*(v),$$
due to the parabolic regularity of  $f$ at  $\bar x$ for  $v.$

{\it Case 2.} $v\in \R^n\backslash \mathcal{A}(\bar x,w).$ Then either $v\not \in \partial f(\bar x)$ or $df(\bar x)(w)\not=\langle v, w \rangle.$
Put
$$\upsilon_t(z):=\frac{f(\bar x+tw+\frac{1}{2}t^2z)-f(\bar x)-tdf(\bar x)(w)}{\frac{1}{2}t^2}\quad  \forall  z\in \R^n, t>0.$$
We have
$$\upsilon_t^*(v)=\frac{f(\bar x)+f^*(v)-\langle v, \bar x\rangle}{\frac{1}{2}t^2}+\frac{df(\bar x)(w)-\langle v, w\rangle}{\frac{1}{2}t}\quad \forall  v\in \R^n, t>0.$$
Since $f$ is parabolically epi-differentiable at $\bar x$ for $w,$  by \cite[Example 13.59]{RW98}, $\epi \upsilon_t$ converges to $\epi  d^2f(\bar x)(w|\cdot)$ as $t\downarrow 0.$
Noting that  $\upsilon_t(\cdot)$ and $d^2f(\bar x)(w|\cdot)$ are proper lower semicontinuous and convex functions, by  \cite[Theorem 11.34]{RW98}, the latter implies that  $\epi \upsilon_t^*$ converges to $\epi  d^2f(\bar x)(w|\cdot)^*$ as $t\downarrow 0.$  So, for any sequence $t_k\downarrow 0,$ by  \cite[Proposition 7.2]{RW98}, there exists a sequence $v_k\to v$ such that
$$ d^2f(\bar x)(w|\cdot)^*(v)=\lim\limits_{k\to \infty}\upsilon_{t_k}^*(v_k).$$
If $v\not\in \partial f(\bar x)$ then
$f(\bar x)+f^*(v)-\langle v,\bar x\rangle>0.$
Thus, by lower semicontinuity of $f^*$,   we see that
$$\begin{array}{rl}d^2f(\bar x)(w|\cdot)^*(v)&=\lim\limits_{k\to \infty}\upsilon_{t_k}^*(v_k)\\
&=\lim\limits_{k\to \infty}\frac{2}{t_k}\left( \frac{f(\bar x)+f^*(v_k)-\langle v_k, \bar x\rangle}{t_k}+df(\bar x)(w)-\langle v_k, w\rangle\right)\\
&=\infty.\end{array}$$
If  $df(\bar x)(w)\not=\langle v, w \rangle$  then, by \eqref{eqSsub} and  \cite[Proposition 13.5]{RW98},
$\langle v, w\rangle < df(\bar x)(w).$  On the other hand, we have
$$f(\bar x)+f^*(v_k)-\langle v_k, \bar x\rangle= f(\bar x)+\sup\limits_{x\in \R^n}[\langle v_k, x\rangle -f(x)]-\langle v_k, \bar x\rangle\geq 0\quad \forall k.$$
Therefore,
$$\begin{array}{rl}d^2f(\bar x)(w|\cdot)^*(v)&=\lim\limits_{k\to \infty}\upsilon_{t_k}^*(v_k)\\
&=\lim\limits_{k\to \infty}\left( \frac{f(\bar x)+f^*(v_k)-\langle v_k, \bar x\rangle}{\frac{1}{2}t_k^2}+\frac{df(\bar x)(w)-\langle v_k, w\rangle}{\frac{1}{2}t_k}\right)\\
&\geq  \lim\limits_{k\to \infty}\frac{df(\bar x)(w)-\langle v_k, w\rangle}{\frac{1}{2}t_k}
=\infty.\end{array}$$
So, we arrive at  the desired conclusion.  $\hfill\Box$

Following Mohammadi and Sarabi \cite{MS20}, we say that function $\psi(x):=g\circ F$  with $(\bar x,\bar v)\in \gph \partial \psi$  satisfies the basic assumptions at $(\bar x, \bar v)$ if  the following conditions hold:
\begin{itemize}
\item[\quad(H1)] the metric subregularity  qualification condition \eqref{MSeq4.5a} is valid at $\bar x$;
\item[\quad(H2)] for each $y\in \Lambda(\bar x,\bar v),$ $g$ is parabolically epi-differentiable at $F(\bar x)$ for every $u\in K_g\big(F(\bar x), y\big);$
\item[\quad(H3)]  $g$ is parabolically regular at $F(\bar x)$ for every  $y\in \Lambda(\bar x,\bar v).$
\end{itemize}
Here $$\Lambda(\bar x,\bar v):=\big\{y\in \partial g\big(F(\bar x)\big) \ | \ \nabla F(\bar x)^*y=\bar v\big\},$$  and
$$K_g\big( F(\bar x), y\big):=\big\{w\in \R^m\ |\ dg\big( F(\bar x)\big)(w)=\langle \bar v, w\rangle\big\}$$ are  the set of {\it Lagrangian multipliers} associated with $(\bar x,\bar v),$ and  the {\it critical cone} of $g$ at  $( F(\bar x), y\big),$ respectively.


Let us consider the following optimization problem:
\begin{equation}\label{MSeq5.5}
\min\limits_{x\in \R^n}-\langle z, \bar v\rangle+d^2g\big(F(\bar x)\big)\big(\nabla F(\bar x)w|\nabla F(\bar x)z+\nabla^2F(\bar x)(w,w)\big).
\end{equation}

\begin{proposition}\label{MSthm5.2} Let $\psi:\R^n\to \overline{\R}$ be a function that is  represented as  \eqref{eq4.3} with the composition $\psi=g\circ F$ satisfying the basic assumptions $(H1)$-$(H3)$ at $(\bar x,\bar v).$  Then the following assertions hold:
\par $(i)$  For each  $w\in K_\psi(\bar x,\bar v),$ the dual problem of   \eqref{MSeq5.5} is
\begin{equation}\label{MSeq5.6}
\max\limits_{y\in \Lambda(\bar x,\bar v)}\left\langle y, \nabla^2F(\bar x)(w,w)\right\rangle+d^2g\big(F(\bar x),y\big)\big(\nabla F(\bar x)\big);
\end{equation}
 the optimal values of  the primal and dual optimization problems \eqref{MSeq5.5} and \eqref{MSeq5.6} are equal and finite. Furthermore, $\Lambda(\bar x,\bar v, w)\cap \tau\B\not=\emptyset,$ where $\Lambda(\bar x,\bar v, w)$ is the optimal solution set of  \eqref{MSeq5.6} and
\begin{equation}\label{MSeq5.7}\tau:=\kappa\ell\|\nabla F(\bar x)\|+\kappa\|\bar v\|+\ell\end{equation}
with $\ell$ and $\kappa$ given in \eqref{eq4.3} and \eqref{MSeq4.5a}, respectively.
\par$(ii)$ $\psi$ is parabolically regular at $\bar x$ for $\bar v,$ and
\begin{equation}\label{MSeq5.13}
\begin{array}{rl}
d^2\psi(\bar x,\bar v)(w)&=\max\limits_{y\in \Lambda(\bar x,\bar v)}\left\{\big\langle y, \nabla^2F(\bar x)(w,w)\big\rangle+d^2g\big(F(\bar x),y\big)\big(\nabla F(\bar x)w\big)\right\}\\
&=\max\limits_{y\in \Lambda(\bar x,\bar v)\cap (\tau\B)}\left\{\big\langle y, \nabla^2F(\bar x)(w,w)\big\rangle+d^2g\big(F(\bar x),y\big)\big(\nabla F(\bar x)w\big)\right\},
\end{array}
\end{equation}
for every $w\in \R^n,$
where $\tau$ is given by \eqref{MSeq5.7}.
\par$(iii)$ $\psi$ is   twice epi-differentiable  at $\bar x$ for $\bar v.$
\end{proposition}
\noindent{\it Proof.}  Take any $w\in K_\psi(\bar x, \bar v)$ and $y\in \Lambda(\bar x,\bar v).$  Then $d\psi(\bar x)(w)=\la \bar v, w\ra,$ $y\in \partial g\big(F(\bar x)\big)$  and $\nabla F(\bar x)^*y=\bar v.$ So, by Proposition \ref{prop3.6},
$$dg\big(F(\bar x)\big)\big(\nabla F(\bar x)w\big)= \la \bar v, w\ra=\la \nabla F(\bar x)^*y, w\ra=\la y,\nabla F(\bar x)w\ra.$$
This means $\nabla F(\bar x)w\in K_g\big(F(\bar x),y\big).$ By $(H2),$ $g$  is parabolically epi-differentiable at $F(\bar x)$ for  $\nabla F(\bar x)w.$
 Thus, by Lemma \ref{Prop4.6}, the function $d^2g\big(F(\bar x)\big)(\nabla F(\bar x)w|\cdot)$ is a proper lower semicontinuous convex function.
 Hence, from \cite[Example 11.41]{RW98} it follows    that the Fenchel dual problem of \eqref{MSeq5.5} is
 \begin{equation}\label{MSeq5.6A}
\max\limits_{\nabla F(\bar x)^*y=\bar v}\left\langle y, \nabla^2F(\bar x)(w,w)\right\rangle-d^2g\big(F(\bar x)\big)\big(\nabla F(\bar x)w|\cdot)^*(y).
\end{equation}
Pick any $y\in \R^m$ with $\nabla F(\bar x)^*y=\bar v.$
If $y\not\in \partial g\big(F(\bar x)\big)$ then, by Lemma \ref{Prop4.6},
\begin{equation}\label{FConjugate1}d^2g\big(F(\bar x)\big)\big(\nabla F(\bar x)w|\cdot)^*(y)=\infty.\end{equation}
 Otherwise, we get
  $y\in \Lambda(\bar x,\bar v).$ Then, by $(H3),$  $g$ is parabolically regular at $F(\bar x)$ for   $y.$
Note that $y\in \partial g\big(F(\bar x)\big)$ and $dg\big(F(\bar x)\big)\big(\nabla F(\bar x)w\big)=\langle y, \nabla F(\bar x)w\rangle.$   So, by Lemma \ref{Prop4.6}, we see that
\begin{equation}\label{FConjugate2} d^2g\big(F(\bar x)\big)\big(\nabla F(\bar x)w|\cdot)^*(y)=-d^2g\big(F(\bar x),y)(w).\end{equation}
From \eqref{FConjugate1} and \eqref{FConjugate2} it follows that problem \eqref{MSeq5.6A} can be written as problem
\eqref{MSeq5.6}. The rest of the proof $(i)$ runs as the one of  \cite[Theorem 5.2]{MS20}, and  the proof of $(ii)$ is similar to the proof of \cite[Theorem 5.4]{MS20}. So, they are omitted.
Finally, we see that, by~$(ii),$  $\psi$ is  parabolically regular at $\bar x$ for $\bar v\in \partial\psi(\bar x)=\partial_p\psi(\bar x),$ and,  by \cite[Theorem 4.4]{MS20}, $\psi$  is  parabolically epi-differentiable at $\bar x$ for every $w\in K_\psi(\bar x,\bar v).$ Therefore, by \cite[Theorem 3.8]{MS20},  $\psi$ is   twice epi-differentiable  at $\bar x$ for $\bar v.$
$\hfill\Box$
\begin{remark}{\rm Under the assumption of Proposition \ref{MSthm5.2},   $g$ is   parabolically regular at $F(\bar x)$ only  for  all  $y\in \partial g\big(F(\bar x)\big)$ with $\nabla F(\bar x)^*y=\bar v.$ Thus, we cannot apply  \cite[Proposition 4.6]{MS20}  to transforming \eqref{MSeq5.6A} into  \eqref{MSeq5.6},  since  \cite[Proposition 4.6]{MS20} requires the parabolic regularity of $g$ at $F(\bar x)$ for every  $y\in \partial g\big(F(\bar x)\big).$ That is why Lemma  \ref{Prop4.6} is utilized. We note  that the results in Propositions~\ref{prop3.6}\&\ref{MSthm5.2} were established in \cite{MS20} for the case where $F$ is twice differentiable in the classical sense.

}\end{remark}

 \section{Quadratic growth and  strong metric subregularity of the subdifferential}

  Let $f\colon\R^n\to\overline{R}$ and     $\ox\in \dom f.$ We say that   $\bar x$ is  a {\it  strong local minimizer} of~$f$  with modulus $\kk>0$ if there is a number $\gg>0$ such that the following  {\it quadratic growth condition}  (QGC) is satisfied:
 \begin{equation}\label{GC}
f(x)-f(\ox)\ge\frac{\kk}{2}\|x-\ox\|^2\quad \mbox{for all}\quad x\in \B_\gg(\ox).
\end{equation}
 The {\it exact modulus} for QGC of  $f$ at $\ox$ is given  by
\[
{\rm QG}\, (f;\ox):=\sup\big\{\kk>0\;|\; \ox \mbox{ is a strong local minimizer of $f$ with modulus $\kk$}\big\}.
\]

 \begin{lemma}{\rm (\cite[Corollary~3.3]{DMN14}).}\label{DMN14} Let $f\colon \R^n\to\oR$ be a  proper lower semicontinuous function and   let  $\ox\in\dom f$ with  $0\in \partial f(\ox)$. Suppose that  the subgradient mapping $\partial f$ is strongly metrically subregular at $\ox$ for $0$ with modulus $\kk>0$ and there are real numbers $r\in(0,\kk^{-1})$ and $\delta>0$ such that
\begin{eqnarray}\label{3.8}
f(x)\ge f(\ox)-\frac{r}{2}\|x-\ox\|^2\quad\mbox{for all}\quad x\in\B_\delta(\ox).
\end{eqnarray}
Then for any $\al\in (0,\kk^{-1})$, there exists a real number $\eta>0$ such that
\begin{eqnarray}\label{3.9}
f(x)\ge f(\ox)+\frac{\al}{2}\|x-\ox\|^2\quad\mbox{for all}\quad x\in\B_\eta(\ox).
\end{eqnarray}
\end{lemma}

\begin{lemma}{\rm (\cite[Lemma 3.6]{CHNT21}).} \label{lem} Let $h:\R^n \to  \oR$ be  a proper function.  Suppose that $h$ is {\em positively homogenenous of degree} $2$ in the sense that  $h(\lm w)=\lm^2h(w)$ for all $\lm>0$ and $w\in \dom h$.  Then for any $w\in \dom h$ and $z\in \partial h(w)$, we have $\la z,w\ra=2h(w)$.
\end{lemma}
 The following result provides some characterizations of the quadratic growth and the strong metric subregularity of the subdifferential.
 \begin{theorem}\label{thm2} Let $f:\R^n\to\oR$ be the function defined by  $f(x)=\varphi(x)+\psi(x)$ for every $x\in \R^n,$  where
  $\varphi:\R^{n}\to\overline\R$  is  twice differentiable at $\ox$ in the extended sense, $0\in \nabla \varphi(\bar x)+\partial \psi(\bar x),$  and  $\psi:\R^{n}\to\overline\R$  is   subdifferentially continuous, prox-regular, and twice epi-differentiable  at $\ox$ for $-\nabla \varphi(\bar x)$.
     Then the following assertions are equivalent:
	\par	$(i)$ The  quadratic growth condition  \eqref{GC} is satisfied.
	\par $(ii)$ The subgradient mapping $\partial f$ is strongly metrically subregular at $(\ox,0),$ and
\begin{equation}\label{thmeq2}
	\la \nabla^{2}\varphi(\ox)w,w\ra+d^{2}\psi\big(\ox|-\nabla \varphi(\ox)\big)(w)\geq 0 \ \mbox{for all} \ w\in \R^n.
	\end{equation}
	\par $(iii)$ The subgradient mapping $\partial f$ is strongly metrically subregular at $(\ox,0),$ and  $\ox$ is a local minimizer for $f.$
	\par 	$(iv$ For all  $w\in {\rm dom}D\partial \psi\big(\ox| -\nabla \varphi(\bar x)\big)\backslash\{0\}$ and $z\in D\partial \psi\big(\ox| -\nabla \varphi(\bar x)\big)(w),$ we have
	\begin{equation}\label{thmeq1}
	\la \nabla^{2}\varphi(\ox)w,w\ra+\la z,w\ra>0.
	\end{equation}
 \par $(v)$  There exists a real number $c>0$ such that
\begin{equation}\label{SO}\la \nabla^{2}\varphi(\ox)w,w\ra+\la z,w\ra\geq c\|w\|^2,
\end{equation}
for all $w\in {\rm dom}D\partial \psi\big(\ox| -\nabla \varphi(\bar x)\big)$ and $z\in D\partial \psi\big(\ox| -\nabla \varphi(\bar x)\big)(w).$
\par		$(vi)$ For every $w\in \R^n\backslash\{0\},$  we have
	\begin{equation}\label{thmeq2s}
	\la \nabla^{2}\varphi(\ox)w,w\ra+d^{2}\psi\big(\ox|-\nabla \varphi(\ox)\big)(w)>0.
	\end{equation}
If one of the above assertions holds then
\begin{equation}\label{eq}
{\rm QG}(f;\ox)=\inf\left\{\frac{\la \nabla^{2}\varphi(\ox)w,w\ra+\la z,w\ra}{\|w\|^2}\, \Big|\begin{array}{rl}&w\in {\rm dom}D\partial \psi\big(\ox| -\nabla \varphi(\bar x)\big),\\
 &z\in D\partial \psi\big(\ox| -\nabla \varphi(\bar x)\big)(w)\end{array}\right\},
\end{equation}
with the convention that $0/0=\infty$.
\end{theorem}
\noindent{\bf Proof.} Under our assumption, by Theorem~\ref{thm1}, we have
\begin{equation}\label{eq1}
	D\partial(\varphi+\psi)(\ox| 0)(w)=\nabla^{2}\varphi(\ox)(w)+D\partial \psi\big(\ox|-\nabla \varphi(\ox)\big)(w),
  	\end{equation}
and	\begin{equation}\label{eq2}d^{2}(\varphi+\psi)\big(\ox| 0\big)(w)=\big\la w,\nabla^{2}\varphi(\ox)w\big\ra+d^{2}\psi\big(\ox| -\nabla \varphi(\ox)\big)(w),\end{equation}
for every $w\in\R^{n}$.
By \eqref{eq2} and \cite[Theorem 13.24]{RW98}, we see that $(i)\Leftrightarrow (vi)$ and $(iii)\Rightarrow (ii).$

We next prove that $0\in \partial_pf(\bar x).$  Since $\psi$  is   subdifferentially continuous and  prox-regular at $\ox$ for $-\nabla \varphi(\bar x)$, we get
$$\liminf\limits_{x\to \bar x}\frac{\psi(x)-\psi(\bar x)+\langle \nabla \varphi(\bar x), x-\bar x\rangle}{\|x-\bar x\|^2}>-\infty.$$
On the other hand, from the extended  twice differentiability of $\varphi$  at $\ox$,  by \eqref{QE}, it follows that
$$\varphi(x)=\varphi(\bar x)+\langle \nabla \varphi(\bar x), x-\bar x\rangle +\frac{1}{2}\langle x-\bar x, \nabla^2\varphi(\bar x)(x-\bar x)\rangle+o(\|x-\bar x\|^2),$$
which  gives us the following estimations
$$\begin{array}{rl}\liminf\limits_{x\to \bar x}\frac{\varphi(x)-\varphi(\bar x)-\langle \nabla \varphi(\bar x), x-\bar x\rangle}{\|x-\bar x\|^2}&=\liminf\limits_{x\to \bar x}\frac{\frac{1}{2}\langle x-\bar x, \nabla^2\varphi(\bar x)(x-\bar x)\rangle+o(\|x-\bar x\|^2)}{\|x-\bar x\|^2}\\
&\geq -\frac{1}{2}\|\nabla^2\varphi(\bar x)\|>-\infty.\end{array}$$
Therefore,
$$\begin{array}{rl}\liminf\limits_{x\to \bar x}\frac{f(x)-f(\bar x)}{\|x-\bar x\|^2}&=\liminf\limits_{x\to \bar x}\left[\frac{\varphi(x)-\varphi(\bar x)-\langle \nabla \varphi(\bar x), x-\bar x\rangle}{\|x-\bar x\|^2}+\frac{\psi(x)-\psi(\bar x)+\langle \nabla \varphi(\bar x), x-\bar x\rangle}{\|x-\bar x\|^2}\right]\\
&=\liminf\limits_{x\to \bar x}\frac{\varphi(x)-\varphi(\bar x)-\langle \nabla \varphi(\bar x), x-\bar x\rangle}{\|x-\bar x\|^2}+\liminf\limits_{x\to \bar x}\frac{\psi(x)-\psi(\bar x)+\langle \nabla \varphi(\bar x), x-\bar x\rangle}{\|x-\bar x\|^2}\\
&>-\infty.
\end{array}$$
This shows that $0\in \partial_pf(\bar x).$

Hence, by \eqref{eq1} and \cite[Theorem 3.2]{CHNT21}, implication $(v)\Rightarrow(iv)\Rightarrow(iii)$  holds, and
\begin{equation}\label{eqI}
{\rm QG}(f;\ox)\geq\inf\left\{\frac{\la \nabla^{2}\varphi(\ox)w,w\ra+\la z,w\ra}{\|w\|^2}\, \Big|\begin{array}{rl}&w\in {\rm dom}D\partial \psi\big(\ox| -\nabla \varphi(\bar x)\big),\\
 &z\in D\partial \psi\big(\ox| -\nabla \varphi(\bar x)\big)(w)\end{array}\right\}.
\end{equation}

We now prove $(ii)\Rightarrow(i).$ Suppose  $\partial f$ is strongly metrically subregular at $\ox$ for $0$ with modulus $\kk>0,$ and \eqref{thmeq2} holds.
 By \eqref{eq2}, we get
 \begin{equation}\label{thmeq2A}
	d^{2}f\big(\ox|0\big)(w)\geq 0 \,\  \mbox{for all} \ w\in \R^n.
	\end{equation}
Fix an arbitrary $r\in (0,\kk^{-1}).$ Then there exists a real number $\delta>0$ such that
\begin{eqnarray}\label{3.8A}
f(x)\ge f(\ox)-\frac{r}{2}\|x-\ox\|^2\quad\mbox{for all}\quad x\in\B_\delta(\ox).
\end{eqnarray}
Indeed, suppose by contrary that this claim does not hold. Then, for each $k\in \N,$  there exists $x_k\in \B_{1/k}(\ox)$ with
$$f(x_k)< f(\ox)-\frac{r}{2}\|x_k-\ox\|^2.$$
Put $t_k:=\|x_k-\ox\|$ and $w_k:=t^{-1}_k(x_k-\ox)$ for $k\in \N.$  We see that $\_k\downarrow0$ as $k\to \infty.$ Furthermore, passing to a subsequence if necessary, we may assume that $\{w_k\}$ converges to some  $\bar w\in \R^n$ as $k\to \infty.$ So we have
$$\begin{array}{rl}d^{2}f\big(\ox|0\big)(\bar w)&=\liminf\limits_{ \begin{subarray}\quad \,\ t\downarrow 0\\
w\longrightarrow \bar w\end{subarray}}\frac{f(\bar x+t w)-f(\bar x)-\tau \langle 0, w\rangle}{ \frac{1}{2}t^2}\\
&\leq \liminf\limits_{k\to\infty}\frac{f(\bar x+t_k w_k)-f(\bar x)}{ \frac{1}{2}t_k^2}\\
&=\liminf\limits_{k\to\infty}\frac{f(x_k)-f(\bar x)}{ \frac{1}{2}\|x_k-\bar x\|^2}\leq -\frac{r}{2}<0.\end{array}$$
This contradicts \eqref{thmeq2A}. Therefore, there exists a real number $\delta>0$ such that
\eqref{3.8A} holds. By Lemma \ref{DMN14}, the quadratic growth condition  \eqref{GC} holds, and we have $(ii)\Rightarrow(i).$

Finally, we prove $(i)\Rightarrow(v)$ and
\begin{equation}\label{eqII}
{\rm QG}(f;\ox)\leq\inf\left\{\frac{\la \nabla^{2}\varphi(\ox)w,w\ra+\la z,w\ra}{\|w\|^2}\, \Big|\begin{array}{rl}&w\in {\rm dom}D\partial \psi\big(\ox| -\nabla \varphi(\bar x)\big),\\
 &z\in D\partial \psi\big(\ox| -\nabla \varphi(\bar x)\big)(w)\end{array}\right\}.
\end{equation}
Suppose that $\bar x$ is a strong local minimizer with modulus~$\kk$ as in \eqref{GC}.  We derive from \eqref{GC} and  \eqref{ssd} that
\begin{equation}\label{eqH0} d^2f(\bar x| 0)(w)\geq \kappa \|w\|^2\quad \mbox{for all}\ w\in \R^n.\end{equation}
Since $\psi$ is  subdifferentially continuous, prox-regular, and twice epi-differentiable  at $\bar x$ for $-\nabla\varphi(\bar x)\in \partial \psi(\bar x),$ it follows from \eqref{Dh}  that
\begin{equation}\label{eqH1}
D(\partial \psi)\big(\bar x|-\nabla\varphi(\bar x)\big)=\partial h\quad \mbox{with}\quad  h(\cdot):=\dfrac{1}{2} d^2 \psi\big(\bar x|-\nabla\varphi(\bar x)\big)(\cdot)
\end{equation}
Note from \eqref{ssd} and \eqref{eqH0}  that $h$ is  proper and   positively homogenenous of degree $2$. By Lemma~\ref{lem}, for any $z\in D(\partial  \psi)\big(\bar x|-\nabla\varphi(\bar x)\big)(w)=\partial h(w),$  we obtain from \eqref{eqH0} and   \eqref{eqH1} that
\begin{equation}\label{eqH2}
\la z,w\ra=2h(w)=d^2 \psi\big(\bar x|-\nabla\varphi(\bar x)\big)(w).
\end{equation}
Therefore,  for every $w\in {\rm dom}D\partial \psi\big(\ox| -\nabla \varphi(\bar x)\big)$ and $z\in D\partial \psi\big(\ox| -\nabla \varphi(\bar x)\big)(w),$ by \eqref{eq1}, \eqref{eq2}, \eqref{eqH0}, and \eqref{eqH2}, we get
$$\la \nabla^{2}\varphi(\ox)w,w\ra+\la z,w\ra= \la \nabla^{2}\varphi(\ox)w,w\ra+d^2 \psi\big(\bar x|-\nabla\varphi(\bar x)\big)(w) =d^2f(\bar x| 0)(w)\geq \kappa \|w\|^2,$$
which clearly verifies $(v)$ and $$ \kk\leq \inf\left\{\frac{\la \nabla^{2}\varphi(\ox)w,w\ra+\la z,w\ra}{\|w\|^2}\, \Big|\begin{array}{rl}&w\in {\rm dom}D\partial \psi\big(\ox| -\nabla \varphi(\bar x)\big),\\
 &z\in D\partial \psi\big(\ox| -\nabla \varphi(\bar x)\big)(w)\end{array}\right\}.$$
 Since  $\kk$ is an arbitrary modulus of the strong local minimizer $\bar x,$  the latter implies  that \eqref{eqII} holds.
So by \eqref{eqI} and \eqref{eqII} we get \eqref{eq}.  $\hfill\Box$
\begin{remark}{\rm By choosing  $\varphi:=0$, we can get \cite[Theorem 3.7]{CHNT21} from Theorem \ref{thm2}.  In the case where $\varphi$ is twice continuously differentiable and $\psi$ is twice epi-differentiable and convex, other characterizations of the quadratic growth as well as the strong metric subregularity of the subdifferential can be found in \cite[Theorem 7.8]{OM21}.
}\end{remark}

We next consider  the composite optimization problem
\begin{equation}\label{MSeq1.1}
\min\limits_{x\in \R^n}  f(x):=\varphi(x)+g\big(F(x)\big),
\end{equation}
where $\varphi:\R^{n}\to\overline\R$  is  twice differentiable at $\ox$ in the extended sense, $F: \R^n\to \R^m$ is twice differentiable, and $g: \R^m\to \overline{\R}:=(-\infty, +\infty]$ is a proper  lower semicontinuous convex function  Lipschitz continuous  around  $F(\bar x)$ relative to its domain with constaint $\ell\in \R_+.$

The {\it Lagrangian} associated with \eqref{MSeq1.1} is defined by
$$L(x,y)=\varphi(x)+\langle F(x), y\rangle-g^*(y),$$
where $g^*(y):=\sup\limits_{v\in \R^m}[\langle y, v\rangle - g(v)]$ is the Fenchel conjugate of $g$ (see \cite{MS20}).

\begin{corollary}\label{hq1} Let  $0\in \nabla \varphi(\bar x)+\partial \psi(\bar x),$
where $\varphi:\R^{n}\to\overline\R$  is  twice differentiable at $\ox$ in the extended sense, and $\psi:=g\circ F$ with $F: \R^n\to \R^m$ being twice  differentiable at  $\bar x$ in the extended sense and $g: \R^m\to \overline{\R}$ being  a proper  lower semicontinuous convex function  Lipschitz continuous  around  $F(\bar x)$ relative to its domain. Assume that the basic assumptions $(H1)$-$(H3)$ hold for $\psi$ at $(\bar x,\bar v)$ with $\bar v:=-\nabla\varphi(\bar x),$ and $\psi$ is prox-regular at $\bar x$ for $\bar v.$ Then, the following assertions are equivalent:
	\par	$(i)$ The  quadratic growth condition  \eqref{GC} is satisfied.
	\par $(ii)$  $\partial f$ is strongly metrically subregular at $(\ox,0),$ and
\begin{equation*}
	\max\limits_{y\in \Lambda(\bar x,\bar v)}\left\{\big\la\nabla^{2}_{xx}L(\bar x,\bar y)w,w\big\ra+d^{2}g\big(F(\bar x),y\big)\big(\nabla F(\bar x)w\big)\right\}\geq 0	\end{equation*}
for all  $w\in K_\psi(\bar x,\bar v);$
	\par $(iii)$  $\partial f$ is strongly metrically subregular at $(\ox,0),$ and  $\ox$ is a local minimizer of  $f.$
	\par 	$(iv$ For all  $w\in {\rm dom}D\partial \psi\big(\ox| -\nabla \varphi(\bar x)\big)\backslash\{0\}$ and $z\in D\partial \psi\big(\ox| -\nabla \varphi(\bar x)\big)(w),$ we have
	\begin{equation*}
	\la \nabla^{2}\varphi(\ox)w,w\ra+\la z,w\ra>0.
	\end{equation*}
 \par $(v)$  There exists a real number $c>0$ such that
\begin{equation}\label{SO}\la \nabla^{2}\varphi(\ox)w,w\ra+\la z,w\ra\geq c\|w\|^2,
\end{equation}
for all $w\in {\rm dom}D\partial \psi\big(\ox| -\nabla \varphi(\bar x)\big)$ and $z\in D\partial \psi\big(\ox| -\nabla \varphi(\bar x)\big)(w).$
\par		$(vi)$ For every $w\in K_\psi(\bar x,\bar v)\backslash\{0\},$  we have
	\begin{equation*}
	\max\limits_{y\in \Lambda(\bar x,\bar v)}\left\{\big\la\nabla^{2}_{xx}L(\bar x,y)w,w\big\ra+d^{2}g\big(F(\bar x),y\big)\big(\nabla F(\bar x)w\big)\right\}> 0.\end{equation*}
If one of the above assertions holds then
\begin{equation*}
{\rm QG}(f;\ox)=\inf\left\{\frac{\la \nabla^{2}\varphi(\ox)w,w\ra+\la z,w\ra}{\|w\|^2}\, \Big|\begin{array}{rl}&w\in {\rm dom}D\partial \psi\big(\ox| -\nabla \varphi(\bar x)\big),\\
 &z\in D\partial \psi\big(\ox| -\nabla \varphi(\bar x)\big)(w)\end{array}\right\},
\end{equation*}
with the convention that $0/0=\infty$.
\end{corollary}
\noindent{\bf Proof.} Under the given assumption, $\psi$ is prox-regular and subdifferentially continuous  at $\bar x$ for $\bar v,$ and  by Proposition \ref{MSthm5.2}, $\psi$ is twice epi-differentiable at $\bar x$ for $\bar v.$
Furthermore,  since $g$ is Lipschitz continuous relative to its domain and $F$ is Lipschitz continuous around $\bar x,$ the composition  $\psi=g\circ F$ is subdifferentially continuous at $\bar x$ for $\bar v.$  On the other hand, by Proposition \ref{MSthm5.2}, we have
$$d^2\psi(\bar x,\bar v)(w)=\max\limits_{y\in \Lambda(\bar x,\bar v)}\left\{\big\la\nabla^{2}F(\bar x)(w,w)\big\ra+d^{2}g\big(F(\bar x),y\big)\big(\nabla F(\bar x)w\big)\right\} \quad \forall w\in \R^n,$$
which  gives us that
$$ \la \nabla^{2}\varphi(\ox)w,w\ra+d^{2}\psi\big(\ox|\bar v\big)(w)=	\max\limits_{y\in \Lambda(\bar x,\bar v)}\left\{\big\la\nabla^{2}_{xx}L(\bar x,y)w,w\big\ra+d^{2}g\big(F(\bar x),y\big)\big(\nabla F(\bar x)w\big)\right\},$$
for every $w\in \R^n.$ Therefore, noting that  $d^2\psi(\bar x,\bar v)$ is a proper lower semicontinuous function with $\dom d^2\psi(\bar x,\bar v)=K_\psi(\bar x,\bar v),$ we get the desired conclusion by applying Theorem \ref{thm2} to the function $f:=\varphi+\psi$ with $\psi:= g\circ F.$
$\hfill\Box$
\begin{remark}{\rm Under  $(H1)$-$(H3)$,  Mohammadi and Sarabi \cite[Theorem 6.3]{MS20}  showed that $(iii)\Leftrightarrow(vi)$ when   $\varphi$ and $F$ are  twice continuously  differentiable around $\bar x$. Since the latter implies the prox-regularity of  $\psi,$  Corollary \ref{hq1} is   an extension of \cite[Theorem 6.3]{MS20}.
}\end{remark}

\begin{example}{\rm  Consider the following optimization problem:
		\begin{equation}\label{E2NLP} \min\limits_{x\in \R} \varphi(x) +\psi(x),
				\end{equation}
		where $\varphi(x)=2x+g(x)$ with  $g(x)$ being taken from Example \ref{ex}, and $\psi(x):=\delta_{\R^2_-}\circ F(x)$ with $F(x)=\big(F_1(x), F_2(x)\big)$,  $F_1(x)=  -x$ and $F_2(x)=-x^3.$
		By  Example \ref{ex},   $\varphi$ is twice differentiable at $\ox=0$ in the extended sense and not prox-regular at $\ox=0$ for $\bar v=0.$ Put
		$$\Gamma=\{x\in \R\, |\,  F_i(x)\leq 0,\ i=1,2\}=\R_{+}\quad \mbox{and}\quad g(y):=\delta_{\R^2_-}(y).$$
		Then $g$ satisfies $(H2)$ and $(H3).$ Furthermore, we see that
		$$d(x,\dom \psi)=d(x,\Gamma)=\begin{cases}
		0 \ &\mbox{if} \ x\geq0,\\
		\left| x\right| \ &\mbox{if} \ x<0,
		\end{cases}$$
		and
		$$d\big(F(x),\dom g\big)=d\big(F(x),\R^{2}_{-}\big)=\begin{cases}
		0 \ &\mbox{if} \ x\geq0,\\
		\sqrt{x^{2}+x^{6}} \ &\mbox{if} \ x<0,
		\end{cases}$$
		which infers that $d(x,\dom \psi)\leq d\big(F(x),\dom g\big).$   This shows that $(H1)$  holds at $\ox.$\\
		We next prove  that $\bar x$ is a strong local minimizer. Indeed, for all $x\in\Gamma\cap[-1,1]$ and $n\in\mathbb{N}^{*},$ we have
		$$x+x^{10/3}\cos\frac{1}{x}\geq0 \ \mbox{and} \ \frac{(2n+1)(2n^2+2n+1)}{n^3(n+1)^3}x+\frac{1}{(n+1)^3}-\frac{1}{n^3}\geq x^{4}\geq0.$$
		Therefore,  we get
		$$\varphi(x)-\varphi(\ox)\geq x\geq x^{2} \ \mbox{for all} \ x\in\Gamma\cap[-1, 1].$$
		Thus, $\bar x$ is a strong local minimizer. By Corollary~\ref{hq1},  the assertions  $(ii)$-$(vi)$   hold.
	}\end{example}
\section{Conclusion}

We have proved some characterizations of the quadratic growth  and the  strong metric subregularity of the subdifferential of  a function  that  can be represented as the  sum of a function twice differentiable in the extended sense and a  subdifferentially continuous, prox-regular, twice epi-differentiable function. Especially, for such a function, we have shown that   the  quadratic growth, the strong metric subregularity of the subdifferential at a local minimizer,   and the positive definiteness of  the subgradient graphical derivative at a  stationary point are equivalent. Our  results are new even for  the case where the twice differentiability  in the extended sense is replaced by the  twice differentiability  in the classical  sense. In this research direction,  it seems to us  that finding out to which extent   the established   results can be applied to  the analysis of  convergence of  numerical  algorithms  is a very interesting issue \cite{BLN18,HMS21,OM21}, which requires further investigation.
Moreover, in order to widen the range of applications of the obtained results, more researches on the class of functions that are twice differentiable in the extended sense  are needed.


\end{document}